\newcommand{\papertitle}{Fitting a Linear Control Policy to Demonstrations with a Kalman Constraint}

\newif\ifpreprint
\preprinttrue

\ifpreprint
\documentclass[12pt]{article}
\usepackage{fullpage}

\title{\papertitle}
\author{
    Malayandi Palan* \\
    \texttt{malayandi@stanford.edu}
    \and
    Shane Barratt* \\
    \texttt{sbarratt@stanford.edu}
    \and
    Alex McCauley \\
    \texttt{alexmccauley@waymo.com}
    \and
    Dorsa Sadigh \\
    \texttt{dorsa@stanford.edu}
    \and
    Vikas Sindhwani \\
    \texttt{sindhwani@google.com}
    \and
    Stephen P. Boyd \\
    \texttt{boyd@stanford.edu}
}

\newcommand{\citep}{\cite}

\else

\documentclass[]{l4dc2020}

\title[\papertitle{}]{\papertitle}
\usepackage{times}

\author{%
 \Name{Malayandi Palan}\footnotemark[1] \Email{malayandi@stanford.edu}\AND
 \Name{Shane Barratt}\footnotemark[2] \Email{sbarratt@stanford.edu}\AND
 \Name{Stephen P. Boyd}\footnotemark[2] \Email{boyd@stanford.edu}\AND
 \addr \footnotemark[1]Department of Computer Science, Stanford University\\
 \addr \footnotemark[2]Department of Electrical Engineering, Stanford University
}
\fi

\usepackage{graphicx,wrapfig,graphics,psfrag,amsmath,amsfonts,verbatim,xcolor,color,siunitx}

\newcommand{\reals}{{\mbox{\bf R}}}

\newcommand{\symm}{{\mbox{\bf S}}}  

\newcommand{\Tr}{\mathop{\bf Tr}}

\newcommand{\Expect}{\mathop{\bf E{}}}

\newcommand{\argmin}{\mathop{\rm argmin}}

\newcommand{\minimize}{\text{minimize}}
\newcommand{\st}{\text{subject to}}

\newcommand{\eg}{{\it e.g.}}
\newcommand{\ie}{{\it i.e.}}

\newcommand{\BEAS}{\begin{eqnarray*}}
\newcommand{\EEAS}{\end{eqnarray*}}
\newcommand{\BEA}{\begin{eqnarray}}
\newcommand{\EEA}{\end{eqnarray}}
\newcommand{\BEQ}{\begin{equation}}
\newcommand{\EEQ}{\end{equation}}
\newcommand{\BIT}{\begin{itemize}}
\newcommand{\EIT}{\end{itemize}}

\newcommand{\states}{\reals^n}
\newcommand{\inputs}{\reals^m}
\newcommand{\disturbances}{\reals^n}

\newcounter{algorithmctr}[section]
\renewcommand{\thealgorithmctr}{\thesection.\arabic{algorithmctr}}
\newenvironment{algdesc}%
   {\refstepcounter{algorithmctr}\begin{list}{}{%
       \setlength{\rightmargin}{0.03\linewidth}%
       \setlength{\leftmargin}{0.03\linewidth}}%
       \rmfamily\small
       \item[]{\setlength{\parskip}{0ex}\hrulefill\par%
        \nopagebreak{\bfseries\textsf{Algorithm \thealgorithmctr~}}}}%
   {{\setlength{\parskip}{-1ex}\nopagebreak\par\hrulefill} \end{list}}

\begin{document}
\maketitle

\begin{abstract}
We consider the problem of learning a linear control policy for a linear
dynamical system, from demonstrations of an expert regulating the
system. The standard approach to this problem is policy
fitting, which fits a linear policy by minimizing a loss function
between the demonstrations and the policy's outputs plus a
regularization function that encodes prior knowledge. Despite its
simplicity, this method fails to learn policies with low or even finite
cost when there are few demonstrations. We propose to add an additional
constraint to policy fitting, that the
policy is the solution to some LQR problem, \ie, optimal in the
stochastic control sense for some choice of quadratic cost. We refer to
this constraint as a Kalman constraint. Policy fitting with a Kalman
constraint requires solving an optimization problem with convex cost and
bilinear constraints. We propose a heuristic method, based on the
alternating direction method of multipliers (ADMM), to approximately
solve this problem. Numerical experiments demonstrate that adding the
Kalman constraint allows us to learn good, \ie, low cost, policies even
when very few data are available.
\end{abstract}

\ifpreprint \else
\fi

\section{Introduction}

\subsection{Fitting a linear policy to demonstrations}
\let\thefootnote\relax\footnotetext{*Equal contribution.}
Typically, we find a control policy for a task as follows. We first
design a cost function that encodes the desired outcomes of the task,
then find a control policy that minimizes that cost function, and
finally we observe or simulate the behavior of this control policy on
the true system. We repeat this process until we are content with
the control policy's performance, either in simulation or in
the real world.

This procedure (optimization-based control) has
been successfully applied to many tasks
\citep{murray2009optimization}. However, for complex tasks, it is
often difficult to find a cost function that precisely captures the
desired task outcomes and can be optimized effectively
\citep{amodei2016faulty, amodei2016concrete}. For example, in autonomous
driving, it is difficult, if not impossible, to construct a cost
function that reliably generates ``comfortable'' driving behavior. For
such tasks, the established procedure mentioned above is very
expensive, time-consuming, and tedious, if it works at all.

Returning to the autonomous driving example, while it may be difficult
to choose a cost function that captures ``comfortable'' driving
behavior, it is relatively straightforward for human operators to
provide demonstrations of such behavior. Similarly, for many other
tasks, it is easier to collect demonstrations of (nearly) optimal
behavior than it is to define a good cost function. This line of
thought has motivated a long line of research on learning from
demonstrations.

Despite this, there has been comparatively little work on learning from
demonstrations in linear systems. This is surprising, since there are
many practical applications of linear systems. Indeed, many systems can
be modeled as linear systems, and we typically find control policies
for nonlinear systems by first approximating these systems as linear
systems. Much progress in control theory has come from studying linear
systems, and we aim to continue that tradition here by considering the
problem of learning a policy from demonstrations on a linear system.

In this paper, we consider the problem of learning a linear policy for a
known stochastic linear system from demonstrations of an expert
regulating the system (\ie, trying to keep the state and input small). The
simplest method to fit a linear policy to demonstrations is (linear)
policy fitting, where we minimize a loss function that measures our fit
to the demonstrations plus a regularization function over linear
policies. Despite its simplicity, policy fitting can lead to unstable
and highly undesirable linear policies when there are few
demonstrations.

Our key insight is the following: since we are trying to learn a policy
for a linear system, the learned policy should be optimal for
\emph{some} quadratic cost function, \ie, some linear quadratic
regulator (LQR) problem. To standard policy fitting we add a
\emph{Kalman constraint}, which requires that the policy be optimal for
some LQR problem. Our name for this constraint refers to the famous
paper by Kalman, \emph{``When is a control system optimal?''}, which
poses the question of determining when a given linear control policy is
LQR optimal for some choice of weights \citep{kalman1964linear}. This
procedure guarantees that the learned policy will retain all the
desirable properties of optimal policies for LQR problems, such as
stability and robustness.
We can think of the Kalman constraint as a very specific form of 
regularization, one that is highly tuned to learning a linear control
policy that is meant to regulate a system.

We formulate policy fitting with a Kalman constraint as a bi-convex
optimization problem, with a convex objective and bi-affine constraints.
From this formulation we derive a heuristic, based on the alternating
direction method of multipliers (ADMM), that can (approximately) solve
this problem by solving a small number of convex optimization problems.
We show through numerical experiments that this method can recover stable,
low-cost policies using very few demonstrations.

\subsection{Related work}

In learning from demonstrations, the goal is to learn a
policy from noisy observations of the optimal policy. This goal is
shared by a few other bodies of work, namely ``inverse optimal control",
``inverse reinforcement learning'', and ``imitation learning''; each of
these topics vary somewhat in the types of problems that they consider,
the methods that they employ, and their notation, although there is
considerable overlap. Below, we discuss some relevant work from each of
these topics, as well as the related idea of incorporating stability in
system identification. (For a more complete survey, see, \eg,
\cite{argall2009survey,hussein2017imitation}.)

\paragraph{Inverse optimal control.} In optimal control, we are given
the cost function and our goal is to find the optimal policy. In inverse
optimal control, we are given the optimal policy and asked to find the
cost function. This topic dates back to Kalman's seminal work in 1964,
where he characterized a sufficient and necessary condition for a linear
policy to be optimal for a given LQR problem \citep{kalman1964linear}.
(Our idea of using a Kalman constraint to regularize the policy learning
procedure is directly inspired by this work.) More recently, it was
shown that we can recover the cost function associated with an optimal
policy for an LQR problem by solving a particular semidefinite program
(SDP) \cite[\S10.6]{boyd1994linear}. Unlike these methods, we do not
assume access to the optimal policy and our focus is not on recovering
the cost function but on learning an effective policy.

\paragraph{Inverse reinforcement learning.} In inverse reinforcement
learning, the goal is to learn a cost function from (noisy)
demonstrations of the optimal policy. We can then find a policy by
optimizing the learned cost function. (Some argue that by learning the
cost function first, we get the added benefit of interpretability
\citep{ng2000algorithms}.) Unlike in inverse optimal control however,
here, we do not assume access to the system dynamics.
Much of the work in this space considers
systems with a finite number of states and inputs
\citep{ng2000algorithms, abbeel2004apprenticeship,
ramachandran2007bayesian, ziebart2008maximum}. More recent work has
extended this work to the continuous state and input space setting by
leveraging advances in deep learning \citep{wulfmeier2015maximum,
finn2016guided}. These methods demonstrate astonishing results at times
but make very few (if any) assumptions and thus typically require large
numbers of demonstrations to produce sensible results. Instead, our
focus here is specifically on known linear systems and the low data regime.

\paragraph{Imitation learning.} In imitation learning, which we refer to
as policy fitting, the goal is to find a policy directly from (noisy)
demonstrations of the optimal policy. Imitation learning (or direct
policy learning or behavior cloning) typically involves learning a
mapping from states to inputs via supervised learning
\citep{pomerleau1989alvinn, sammut1992learning, kuniyoshi1994learning,
hayes1994robot}. However, standard imitation learning methods are prone
to instability when only a few demonstrations are available
\citep{ross2011reduction} --  a fact we confirm empirically in this
work. Much like standard imitation learning methods, we attempt to learn
a policy directly from states to inputs in this work; however, unlike
prior work, we focus specifically on linear systems, leveraging their
structure to add prior knowledge. Indeed, our method can be interpreted
as a stability-regularized imitation learning method for linear systems.

\paragraph{Stable system identification.}
Another related problem is learning dynamics from measurements of a
dynamical system. The standard approach to this problem, system
identification \citep{ljung1999system}, frames this problem as a
regression task. Recent work in this space has explored the idea
of leveraging prior knowledge that the system to be identified is
stable. For example, in \citep{lacy2003subspace,boots2008constraint} the
problem of fitting dynamics matrices, subject to the constraint that the
dynamics are asymptotically stable, is framed as a convex optimization
problem. This work has also been extended to nonlinear systems
\citep{singh2018learning}.

\subsection{Outline}
In \S\ref{s-lq-problem}, we introduce the problem of learning from
demonstrations on a linear system via policy fitting. In
\S\ref{s-policy-learning}, we introduce the Kalman constraint, combine
policy fitting with the Kalman constraint, and give an approximate
solution method. In \S\ref{s-examples}, we illustrate our method on
several numerical examples. In \S\ref{s-extensions}, we describe some
natural extensions and variations. In \S\ref{s-conclusion}, we conclude.

\section{Linear policy fitting}\label{s-lq-problem}

\paragraph{Dynamics.}
We consider a fully-observable linear dynamical system of the form
\BEQ\label{e-dynamics}
    x_{t+1} = A x_t + B u_t + \omega_t, \quad t=0,1,\ldots,
\EEQ
where, at time $t$,
$x_t\in\states$ is the system state, $u_t\in\inputs$
is the control input, and $\omega_t\in\disturbances$ is a (random)
disturbance. The (known) dynamics matrices are $A \in\reals^{n\times n}$
and $B \in\reals^{n\times m}$, and we assume that $(A,B)$ is
controllable. We assume that $\Expect[\omega_t]=0$,
$\Expect[\omega_t\omega_t^T]=W$, and that $\omega_t$ is independent of
the state $x_t$, the control input $u_t$, and $\omega_\tau$ for
$\tau\neq t$.

\paragraph{Policy.}
A \emph{policy} $\pi: \states \to \inputs$ is a function that maps the
current state $x_t$ to the control input $u_t$ that we will apply
to the system,
\BEQ\label{e-policy}
    u_t = \pi(x_t), \quad t=0,1,\ldots.
\EEQ
Equations \eqref{e-dynamics} and \eqref{e-policy} together define the
closed-loop system.
We will consider the case where $\pi$ is linear, or
\BEQ\label{e-linear-policy}
u_t = \pi(x_t) = K x_t, \quad t=0,1,\ldots,
\EEQ
where $K \in \reals^{m \times n}$ in the gain matrix.
The closed-loop dynamics are then
\BEQ\label{e-closed-loop}
x_{t+1} = (A+BK)x_t + \omega _t, \quad t=0,1,\ldots.
\EEQ

\paragraph{``Expert'' demonstrations.}
We consider the case where we observe some ``expert'' regulating the
system \eqref{e-dynamics}, \ie, trying to keep the state $x_t$ and input
$u_t$ small. We receive $N$ demonstrations, $(x^i,u^i)$, $i=1,\ldots,N$.
These state input pairs need not be ordered in time, optimal in any sense, or even
deterministic (\ie, we can have pairs with the same state and 
different inputs).
We do assume, however, that the expert is attempting in good faith to
regulate the system, \ie, keep $x_t$ and $u_t$ small, in some sense.

\paragraph{Linear policy fitting.}
The goal in linear policy fitting is to fit a linear policy $K$ to the
demonstrations. To that end, we consider the average of some loss
function $l:\reals^m\times\reals^m\to\reals$ (convex in its first
argument) across the demonstrations:
\[
L(K) = \sum_{i=1}^N l(Kx^i, u^i).
\]
We also assume that we have a convex regularization function
$r:\reals^{m \times n} \to \reals \cup \{+\infty\}$ that encodes prior
knowledge on $K$. Infinite values of $r$ can be interpreted as
constraints on $K$.

To fit the policy, we solve the problem
\begin{equation}
\begin{array}{ll}
    \mbox{minimize} & L(K) + r(K),\\
\end{array}
\label{e-policy-fitting-problem}
\end{equation}
with variable $K$. 
This optimization
problem is convex, and so can be solved efficiently
\citep{boyd2004convex}. We denote a solution to
\eqref{e-policy-fitting-problem} by $K^\mathrm{pf}$.

The objective function in \eqref{e-policy-fitting-problem} consists of
two parts: the demonstration loss $L$ and the regularization function
$r$. The first term here encourages a good fit to the demonstrations,
and the second term encourages the policy to be simpler or to be 
consistent with prior knowledge.

If the expert is indeed using a linear policy,
\ie, $x^i=K^\mathrm{expert}u^i$, then we
will recover $K^\mathrm{expert}$ using $n$ demonstrations with high
probability, with any reasonable choice of $l$, and without regularization.
However, in general, policy fitting does not perform well
when there are only a few demonstrations. Often, $K^\mathrm{pf}$ does
not even stabilize the system, let alone mimic the expert's policy well
in closed-loop simulation (see \S\ref{s-examples} and
\cite{ross2011reduction}).

\section{Policy fitting with a Kalman constraint}
\label{s-policy-learning}

\subsection{Linear quadratic regulator}
The well-known linear quadratic regulator (LQR) problem
\citep{bertsekas2017dynamic} chooses a policy $\pi$ that minimizes the
average of a quadratic cost function,
\BEQ
J(\pi) = \lim_{T\to\infty}\frac{1}{T}\Expect\sum_{t=0}^{T-1}
\left( x_t^TQx_t + u_t^TRu_t\right),
\label{e-cost}
\EEQ
where $Q \in \symm_{+}^n$ (the set of symmetric $n\times n$ positive
semi-definite matrices) and $R \in \symm_{++}^m$ (the set of symmetric
$m \times m$ positive definite matrices) are the state and control
weight matrices, respectively, subject to the dynamics
\eqref{e-dynamics} and $u_t=\pi(x_t)$. Here the expectation is taken
over the initial state $x_0$ and disturbances $\omega_0,\omega_1,\ldots$.

It is well known that (assuming some reasonable technical conditions hold) 
the optimal policy $\pi$ is linear, of the form
\[
\pi^\star(x_t) = K^\star x_t,
\]
where $K^\star\in\reals^{m \times n}$ is the \emph{optimal gain matrix}.
The optimal gain matrix $K^\star$ depends on $Q$ and $R$, as well as $A$ and $B$,
but not $W$.

\paragraph{When is a linear policy optimal for given $Q$ and $R$?}
Suppose we are given a linear policy $K$, and want to know if $K$ is
optimal for the cost \eqref{e-cost}, given the quadratic cost matrices
$Q$ and $R$. This is the case when
\[
K = -(R+B^TPB)^{-1}B^TPA,
\]
where $P$ is the unique positive semi-definite solution of
the algebraic Riccati equation
\[
P = Q + A^TPA - A^TPB(R+B^TPB)^{-1}B^TPA.
\]

\paragraph{When is a linear policy optimal for \emph{some} $Q$ and $R$?}
Suppose we are given a linear policy $K$, and want to know if $K$ is
optimal for the cost \eqref{e-cost} for \emph{some} quadratic cost
matrices $Q$ and $R$, which we are free to choose. This question was
addressed in \cite[\S10.6]{boyd1994linear}, where the authors showed
that we can answer this question, and get $Q$ and $R$, by solving the
(convex) semidefinite feasibility problem
\begin{equation}
\begin{array}{ll}
    \mbox{minimize} & 0,\\
    \mbox{subject to} & Q + A^TP(A + BK) - P = 0,\\
    & RK + B^TP(A + BK) = 0,\\
    & P \succeq 0, \quad Q \succeq 0, \quad R \succ 0,\\
\end{array}
\label{e-sdp}
\end{equation}
with variables $P$, $Q$, and $R$, where $\succeq$ denotes 
matrix inequality, \ie, with respect to the semidefinite cone.
If problem \eqref{e-sdp} is feasible, then $K$
is optimal for the quadratic cost matrices $Q$ and $R$. On the other
hand, if \eqref{e-sdp} is infeasible, then $K$ is not optimal for any
LQR problem.

\paragraph{Kalman constraint.}
We refer to the constraints in \eqref{e-sdp} as a \emph{Kalman
constraint} on $K$, in tribute to Kalman's seminal work on optimal
control, entitled ``\emph{When is a linear control system optimal?}''
\citep{kalman1964linear}. A Kalman constraint on $K$ implies that $K$
must be optimal for some quadratic cost matrices $Q$ and $R$. Policies
that satisfy a Kalman constraint retain all the desirable properties of
an LQR-optimal policy, such as stability and robustness.

\subsection{Policy fitting with a Kalman constraint}
We add a Kalman constraint to the regularization function $r$ in the
standard policy fitting problem \eqref{e-policy-fitting-problem}. This
new policy fitting problem has the form
\begin{equation}
\begin{array}{ll}
    \mbox{minimize} & L(K) + r(K),\\
    \mbox{subject to} & Q + A^TP(A + BK) - P = 0,\\
    & RK + B^TP(A + BK) = 0,\\
    & P \succeq 0, \quad Q \succeq 0, \quad R \succeq I,\\
\end{array}
\label{e-non-convex}
\end{equation}
with variables $P$, $Q$, $R$, and $K$. Note that the optimal gain matrix
$K^\star$ is invariant to the relative scale of the $(Q, R)$ matrices,
\ie, if $K$ is optimal for $(Q,R)$, it is also optimal for $(\alpha Q,
\alpha R)$ for any $\alpha > 0$. Thus, we can replace the (open)
constraint $R \succ 0$ in \eqref{e-sdp} with the (closed) constraint $R
\succeq I$ by suitable scaling.

This problem is nonconvex and so, in general, difficult to solve
exactly. However, we note that this problem is bi-convex in $(P,Q,R)$ and
$K$. In \S\ref{s-solution}, we derive a heuristic method that finds an
approximate or local solution.

\paragraph{Interpretability.}
Previous work in inverse reinforcement learning suggests that by
recovering the (unknown) cost function, we will be able to better
interpret the expert's policy \citep{ng2000algorithms}. Indeed, our
method does recover a cost function that, at least approximately,
explains the observed demonstrations. However, we note that, even in
linear systems, for a given linear policy, the quadratic cost function
is not unique, even up to a scale factor.
The problem of recovering a cost function from a given policy is
under-determined, so even if we recover the true policy, the stage cost
coefficients $Q$ and $R$, and the cost-to-go matrix $P$ may not converge
to the true $Q$ and $R$, if they even exist. Therefore, all we can hope
to do is recover a stable and desirable policy.
(If we do care about the recovered $Q$ and $R$, and have some
prior knowledge about them, we can add a suitable regularization term
on $Q$, $R$, or $P$.)

\subsection{Solution method}
\label{s-solution}
\paragraph{ADMM.}
We observe that the objective in \eqref{e-non-convex} is convex, and
that the bi-convexity of the problem comes from the bi-affine constraints.
Therefore, we propose to use the alternating direction method of
multipliers (ADMM), which is guaranteed to converge to a (not
necessarily optimal) stationary point for this problem, provided the
penalty parameter is sufficiently large \citep{gao2019admm}.

The augmented Lagrangian \citep{hestenes1969multiplier} of
\eqref{e-non-convex} is the extended function
\[
\mathcal L_\rho(K, P, Q, R, Y) = L(K) + r(K) + I_C(P, Q, R) + \Tr(Y^TM) +
\frac{\rho}{2} \|M\|_F^2,
\]
where $Y$ is the dual variable for the constraints in
\eqref{e-non-convex},
\[
M = \begin{bmatrix}Q + A^TP(A+BK) - P\\ RK + B^TP(BK+A)\end{bmatrix},
\quad I_C(P, Q, R) = \begin{cases} 0 & P \succeq 0, Q \succeq 0, R
\succeq I, \\ +\infty & \text{otherwise},
\end{cases}
\]
and $\rho>0$ is an algorithm parameter.

The ADMM algorithm alternates between minimizing the augmented Lagrangian over $K$,
minimizing the augmented Lagrangian over $(P, Q, R)$, and performing a
dual update to $Y$. The full procedure is summarized in algorithm
\ref{alg:lfd} below.
\begin{algdesc}
\label{alg:lfd}
\emph{Learning from demonstrations in linear systems via ADMM.}
\begin{tabbing}
    {\bf given} initial parameters $K^0, P^0, Q^0, R^0$, penalty
    parameter $\rho$, iterations $n_\mathrm{iter}$.\\
    {\bf for} $k=0,\ldots,n_\mathrm{iter}-1$ \\
        \qquad \=\ 1.\ \emph{$K$ step.} Let $K^{k+1} = \argmin_K
        \mathcal L_{\rho}(K, P^k, Q^k, R^k, Y^k)$.\\[.5em]
        \qquad \=\ 2.\ \emph{$(P, Q, R)$ step.} Let
        $(P^{k+1},Q^{k+1},R^{k+1}) = \argmin_{(P, Q, R)} \mathcal
        L_{\rho}(K^{k+1}, P, Q, R, Y^k)$.\\[.5em]
        \qquad \=\ 3.\ \emph{$Y$ step.} Let $Y^{k+1} = Y^k + \rho
        \begin{bmatrix}Q^{k+1} + A^TP^{k+1}(A+BK^{k+1}) - P^{k+1}\\
        R^{k+1}K^{k+1} + B^TP^{k+1}(BK^{k+1}+A)\end{bmatrix}$.\\
    {\bf end for}\\
    {\bf return} $K^{n_\mathrm{iter}}$.
\end{tabbing}
\end{algdesc}

\paragraph{$K$ step.} The update for $K$ can be expressed as the
solution to the convex optimization problem
\BEA
    \minimize \quad &&L(K) + r(K) +
    \frac{\rho}{2}\left\|\begin{bmatrix}A^TP^kBK - (P^k - Q^k - A^TP^kA
    - \frac{1}{\rho}Y_1^k)\\ (R^k + B^TP^kB) K - (-B^TP^kA -
    \frac{1}{\rho}Y_2^k)\end{bmatrix}\right\|_F^2\label{eq:step1}
\EEA
where $Y^k=(Y_1^k, Y_2^k)$.
This step can be interpreted as performing policy fitting with an
additional term in the regularization function
that suggests that $K$ should be approximately optimal for the LQR
control problem with cost matrices $Q^k$, $R^k$, and
cost-to-go matrix $P^k$.

\paragraph{$(P, Q, R)$ step.} The update for $(P, Q, R)$ can be
expressed as the solution to the convex optimization problem
\BEA
    \minimize \quad &&\left\|\begin{bmatrix}
    Q + A^TP(A+BK^{k+1}) - P + \frac{1}{\rho}Y_1^k\\
    RK^{k+1} + B^TP(BK^{k+1}+A) + \frac{1}{\rho}Y_2^k
    \end{bmatrix}\right\|_F^2\nonumber\\
    \st \quad &&P\succeq 0, Q \succeq 0, R \succeq I,\label{eq:step2}
\EEA
which can be interpreted as finding the cost matrices and
cost-to-go-matrix of an LQR problem such that $K^{k+1}$ is
approximately optimal for that problem.

\paragraph{Termination criterion.} We can either run the algorithm
for a fixed number of iterations ($n_\mathrm{iter}$ in algorithm \ref{alg:lfd})
or until the Frobenius between successive values of $K$ is less than
some chosen value $\epsilon$.

\paragraph{Convergence.} This algorithm is only guaranteed to converge
if the penalty parameter is large enough \citep{gao2019admm};
even when it does converge,
it need not converge to an optimal value. It is simply a sophisticated
heuristic for finding an effective, stable policy.

\section{Examples}\label{s-examples}

We illustrate our method and compare it with linear regression on
several problems. In all of our experiments we use $\rho=1$, and run ADMM
with a zero initialization and 5 random initializations, ultimately using the $K$
with the lowest value of $L(K) + r(K)$.
We use CVXPY~\citep{diamond2016cvxpy} to implement
algorithm \ref{alg:lfd}.

\subsection{Imperfect LQR}
\label{sec:imperfect_lqr}
We first consider the case where the expert
is performing imperfect regulation in an LQR problem.
That is, our demonstrations have the form
\[
u^i = K^\star x^i + z^i, \quad z^i \sim \mathcal N(0, \Sigma), \quad i=1,\ldots,N,
\]
where $K^\star$ is the solution to an LQR problem with
cost matrices $Q^\mathrm{true}$ and $R^\mathrm{true}$
and $\Sigma \succ 0$.
We consider loss and regularization functions
\[
l(\hat u, u) = \|\hat u - u\|_2^2, \quad r(K) = (0.01)\|K\|_F^2.
\]

\begin{figure}
    \centering
    \includegraphics[]{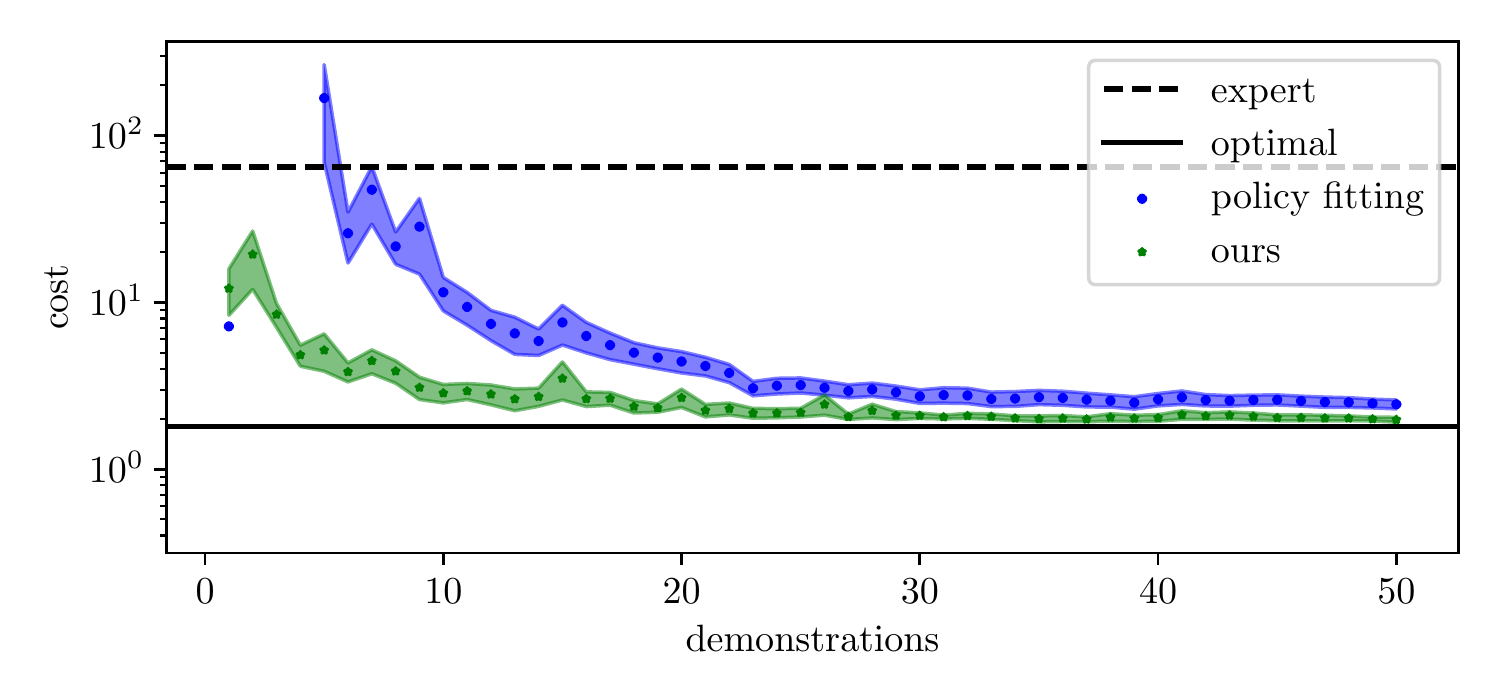}
    \caption{Small random example: The expected cost of policy fitting
    and our method for a range of numbers of demonstrations. Infinite values
    are ignored. }
    \label{fig:small_random}
\end{figure}

\begin{figure}
    \centering
    \includegraphics[]{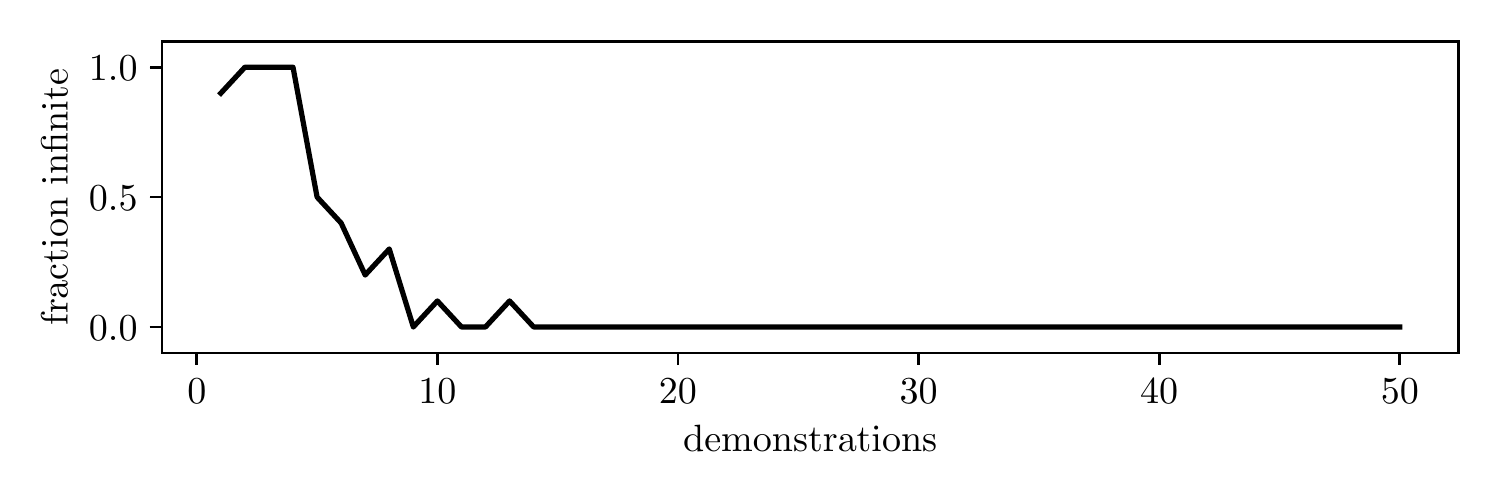}
    \caption{Small random example: The fraction of simulations for which the
    cost for policy fitting was infinite. Our method always attained finite
    cost.}
    \label{fig:inf_small_random}
\end{figure}

\paragraph{Small random example.}
We consider a system with $n=4$ states and $m=2$ inputs. The data is
generated according to
\[
A_{ij} \sim \mathcal N(0, 1), \quad B_{ij} \sim \mathcal N(0, 1), \quad
Q^\mathrm{true}=I, \quad R^\mathrm{true}=I, \quad W=(0.25)I, \quad \Sigma=(4)I,
\]
and the matrix $A$ is scaled so that its spectral radius is one.
We ran policy fitting
with and without a Kalman constraint on varying numbers of demonstrations,
and averaged the results over ten random seeds.
In figure \ref{fig:small_random} we show the expected cost (when it is finite) of
policy fitting and our method versus the number of demonstrations,
as well as the expected cost incurred by the expert and the optimal policy.
Whereas our method never incurred infinite cost,
standard policy fitting did; figure \ref{fig:inf_small_random} shows
the fraction of the time that the cost for policy fitting was finite.


\paragraph{Aircraft example.}
We consider the control of a 747 aircraft during level flight at an elevation of 40000
feet, traveling at 774 feet per second. The states and inputs represent
deviations from operating or trim conditions. The states are $u$, the
velocity of the aircraft along the body axis (in ft/s), $v$, the
velocity of the aircraft perpendicular to the body axis (in ft/s),
$\theta$, the angle between the body axis and horizontal (in crad), and
$q=\dot\theta$, the pitch rate (in crad/s). The inputs are $\delta_e$,
the elevator angle (in crad), and $\delta_t$, the thrust.

\begin{figure}
    \centering
    \includegraphics[]{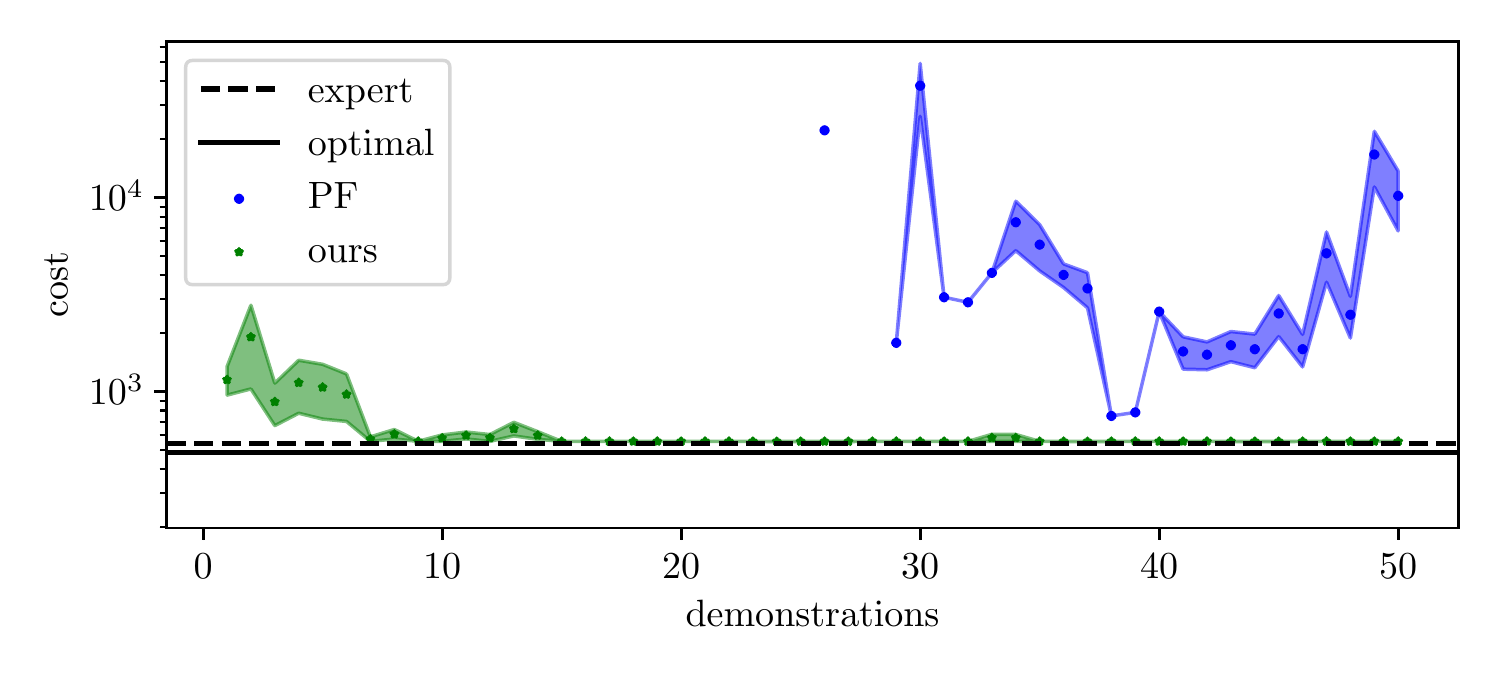}
    \caption{Aircraft control example: The expected cost of policy fitting
    and our method for a range of numbers of demonstrations. Infinite values
    are ignored.}
    \label{fig:aircraft}
\end{figure}

\begin{figure}
    \centering
    \includegraphics[]{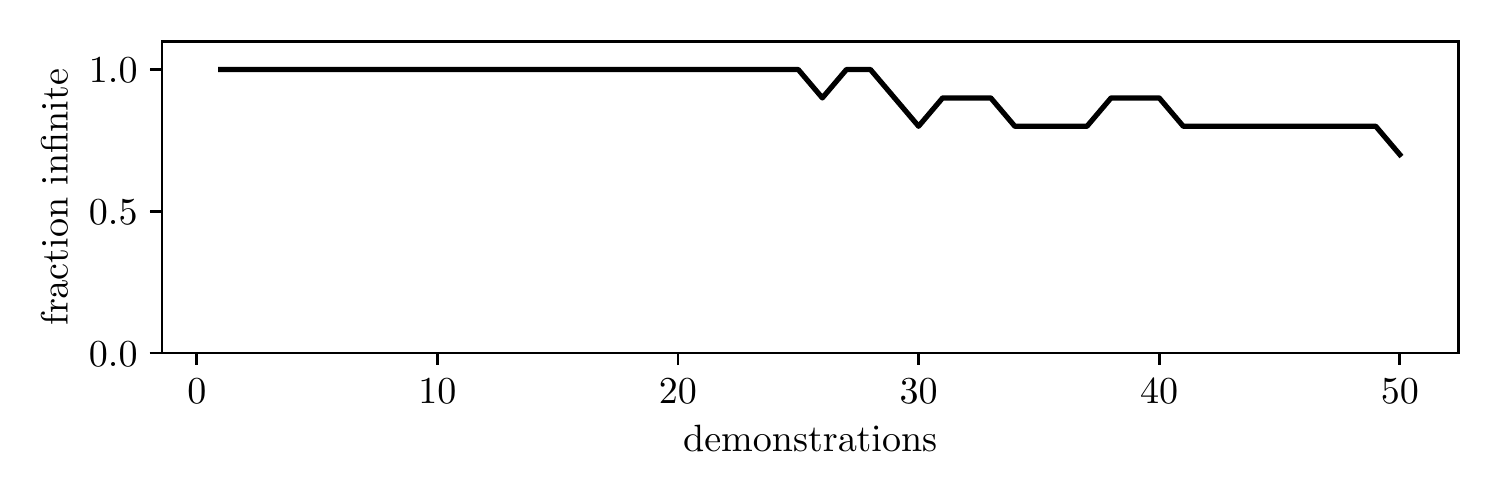}
    \caption{Aircraft control example: The fraction of simulations for which
    the cost for policy fitting was infinite. Our method always attained finite
    cost.}
    \label{fig:inf_aircraft}
\end{figure}

The linearized dynamics, discretized at an interval of 0.01 seconds,
have the form
\[
\begin{bmatrix}
u_{t+1} \\ v_{t+1} \\ q_{t+1} \\ \theta_{t+1}
\end{bmatrix} = 
\begin{bmatrix}
1 & 0.039 & 0 & -0.322 \\
-0.065 & 0.997 & 7.74 & 0 \\
0.02 & -0.101 & 0.996 & 0 \\
0 & 0 & 1 & 1
\end{bmatrix}
\begin{bmatrix}
u_t \\ v_t \\ q_t \\ \theta_t
\end{bmatrix}
+ \begin{bmatrix}
.0001 & 0 \\
-.0018 & -.0004 \\
-.0116 & .00598 \\
0 & 0
\end{bmatrix}
\begin{bmatrix}
(\delta_e)_t \\
(\delta_t)_t
\end{bmatrix} + \omega_t,
\]
where the disturbance $\omega_t$ is caused by wind, with covariance
\[
W = 
\begin{bmatrix}
 0.100 & -0.003 &  0.002   &  0       \\
-0.003 &  0.1 & -0.010 &  0       \\
 0.002   & -0.010 &  0.001   &  0       \\
 0        &  0        &  0        &  0       
\end{bmatrix}.
\]
We use the cost matrices $Q=I$, $R=I$, incentivizing us to
keep the 747 at the trim condition while keeping the thrust and elevator
angle at the nominal levels.
We also use the observation noise $\Sigma=(25)I$.
We ran policy fitting
with and without a Kalman constraint on varying numbers of demonstrations,
and averaged the results over ten random seeds.
In figure \ref{fig:aircraft} we show the expected cost (when it is finite) of
policy fitting (PF) and our method versus the number of demonstrations,
as well as the expected cost incurred by the expert and the optimal policy.
Whereas our method never incurred infinite cost,
standard policy fitting did; figure \ref{fig:inf_aircraft} shows
the fraction of the time that the cost for policy fitting was finite.


\subsection{LQR with outliers}
\label{sec:lqr_outliers}

We consider an LQR problem, with demonstrations generated according to
\[
u^i = K^\star x^i + z^i, \quad z^i \sim \mathcal N(0, \Sigma), \quad i=1,\ldots,N,
\]
where $K^\star$ is the solution to an LQR problem with
cost matrices $Q^\mathrm{true}$ and $R^\mathrm{true}$,
and $\Sigma \succ 0$.
For each entry of $u^i$, $i=1,\ldots,N$, we flip its sign
with probability $0.1$.
(The entries that are flipped are called outliers.)

We employ the following loss and regularization functions
\[
l(\hat u, u) =\sum_i \phi(\hat u_i - u_i), \quad r(K) = (0.01)\|K\|_F^2,
\]
where $\phi$ is the Huber penalty function, with parameter $M$,
\[
\phi(a) = \begin{cases} a^2/2 & |a| \leq M\\
M|u| - M^2/2 & |a| > M.
\end{cases}
\]
(We use the Huber loss function because it is robust to outliers
\citep[\S 6.1]{boyd2004convex}.)

\begin{figure}
    \centering
    \includegraphics[]{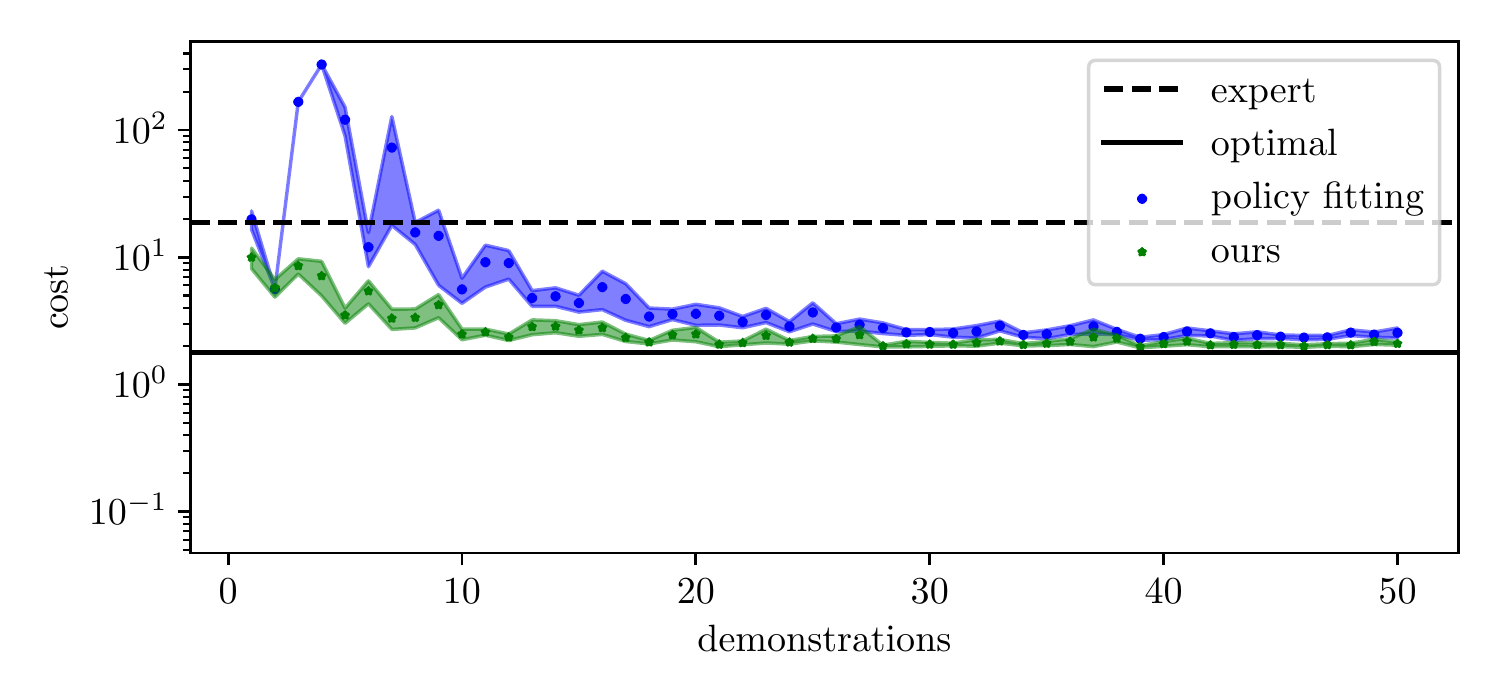}
    \caption{LQR with outliers: The expected cost of policy fitting
    and our method for a range of numbers of demonstrations. Infinite values
    are ignored.}
    \label{fig:lqr_outliers}
\end{figure}

\begin{figure}
    \centering
    \includegraphics[]{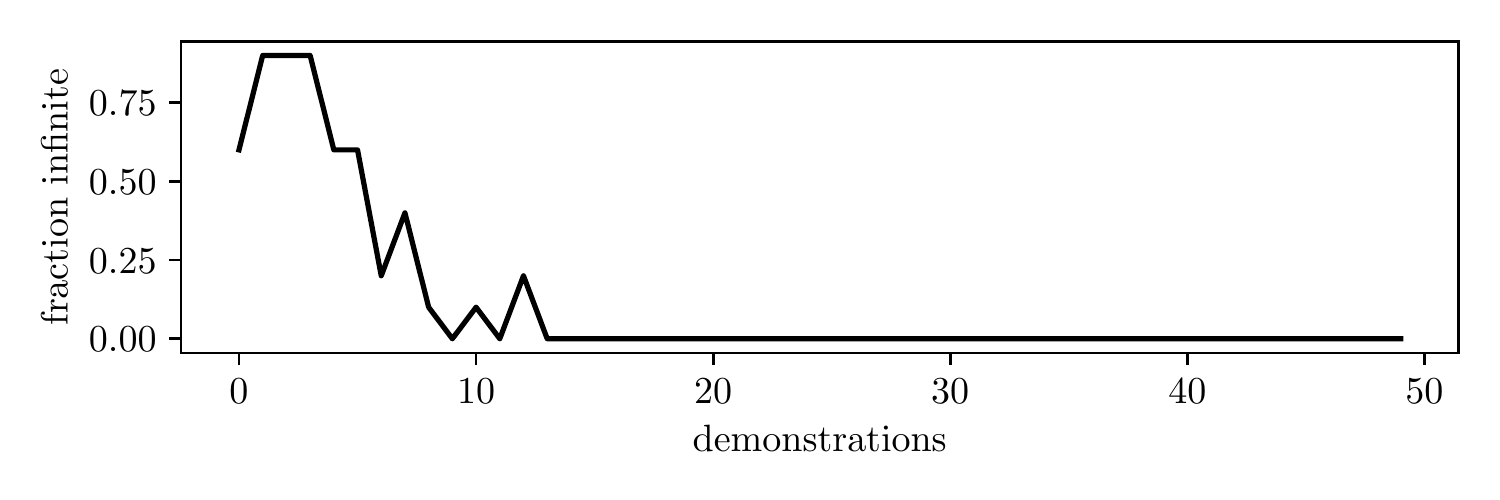}
    \caption{LQR with outliers: The fraction of simulations for which
    the cost for policy fitting was infinite. Our method always attained finite
    cost.}
    \label{fig:inf_lqr_outliers}
\end{figure}

\paragraph{Numerical example.}
We consider the same data as
the small random problem in \S\ref{sec:imperfect_lqr},
and with $M=0.5$.
We ran policy fitting
with and without a Kalman constraint on varying numbers of demonstrations,
and averaged the results over ten random seeds.
In figure \ref{fig:lqr_outliers} we show the expected cost (when it is finite) of
policy fitting and our method versus the number of demonstrations,
as well as the expected cost incurred by the expert and the optimal policy.
Whereas our method never incurred infinite cost,
standard policy fitting did; figure \ref{fig:inf_lqr_outliers} shows
the fraction of the time that the cost for policy fitting was finite.


\section{Extensions and variations}\label{s-extensions}
\paragraph{General quadratic cost problem.} In this paper, in the
interest of clarity, we focused on the regulation of linear systems,
where the goal is to keep the state and input small. However, our method
can be easily adapted for a more general class of problems, such as
tracking problems, where the goal is not to keep the state and input
small but to keep the state close to a given trajectory. To do so, we
simply need to define the quadratic stage cost in its more general form
\[
\begin{bmatrix}x_t\\u_t\\1\end{bmatrix}^T Q \begin{bmatrix}x_t\\u_t\\1\end{bmatrix},
\]
where $Q\in\symm_{+}^{n+m+1}$.
We can then
replace the constraints in
\eqref{e-non-convex} with the appropriate Riccati equation.

\paragraph{Finite-horizon problem.} Similarly, in this paper, we only
considered the time-invariant, infinite horizon problem, where the
dynamics function is as given in \eqref{e-dynamics}. However, our method
can be easily extended for time-varying, finite-horizon problems, where
the dynamics function is given by
\[
x_{t+1} = A_tx_t + B_tu_t + \omega_t, \quad t=0, 1, \ldots, T-1,
\]
where $T$ is the horizon of the problem. (In this problem, the cost
function is similarly truncated to T steps and is also time-varying).
Here, our goal is to learn a policy for each time-step, $K^\star_0,
K^\star_1,\ldots, K^\star_T$ instead of a single policy.

To solve this problem using our method, we first need replace the
constraints in \eqref{e-non-convex} with constraints for each time-step.
(There are thus $T$ times as many constraints.) We can then adapt
our algorithm to solve this problem by following the same steps outlined
in \S\ref{s-solution}.

\section{Conclusion}\label{s-conclusion} In this paper, we introduced a
method for learning policies from demonstrations in linear
systems and showed in numerical experiments that this method outperforms
a widely-used baseline. Our method, which is based on convex
optimization, is easy to implement and consistently produces reliable
results, in contrast to gradient-based methods that are difficult to
make work. We believe that this method and its extensions (see
\S\ref{s-extensions}) have wide-ranging practical applications,
especially in the domain of autonomous driving; indeed, a rigorous
examination of this claim is the subject of future work. We are very
optimistic about the potential of convex optimization to solve modern
control problems and believe that this space will re-emerge as a fruitful
area for research in the coming years.

\section*{Acknowledgments}
S. Barratt is supported by the National Science Foundation Graduate Research Fellowship
under Grant No. DGE-1656518.

\bibliography{lfd_lqr}

\end{document}